\documentclass[11pt,letterpaper]{article}
\usepackage{graphicx} 
\usepackage[margin=1in]{geometry}

\usepackage{amsmath}
\usepackage{amssymb}
\usepackage{amsthm}
\usepackage{hyperref}
\usepackage{color}

\usepackage{enumitem}
\setlist{leftmargin=6mm,nolistsep,noitemsep}

\usepackage[capitalize,noabbrev]{cleveref}
\crefname{question}{Question}{Questions}
\crefname{step}{Step}{Steps}
\crefname{claim}{Claim}{Claims}
\crefname{problem}{Problem}{Problems}
\crefname{observation}{Observation}{Observations}

\newcommand{\MP}{{\rm MP}}

\newcommand{\ie}{i.e., }
\renewcommand{\S}{\mathcal S}

\def\P{\mathcal{P}}
\def\NP{\mathcal{NP}}
\def\BPP{\mathcal{BPP}}

\def\01{\ensuremath{0\mathord{-}1}}
\allowdisplaybreaks

\newtheorem{theorem}{Theorem}

\allowdisplaybreaks

\DeclareMathOperator{\poly}{poly}

\newcommand{\LP}{{\rm LP}}

\usepackage[title]{appendix}

\title{Extended formulations for the multilinear polytope of acyclic hypergraphs}
\author{Alberto Del Pia
\thanks{Department of Industrial and Systems Engineering \& Wisconsin Institute for Discovery,
             University of Wisconsin-Madison.
             E-mail: {\tt delpia@wisc.edu}.
             }
\and
Aida Khajavirad
\thanks{Department of Industrial and Systems Engineering,
             Lehigh University.
             E-mail: {\tt aida@lehigh.edu}.
             }
}
\date{January 8, 2025}

\begin{document}

\maketitle

\begin{abstract}
This article provides an overview of our joint work on binary polynomial optimization over the past decade. We define the multilinear polytope as the convex hull of the feasible region of a linearized binary polynomial optimization problem. By representing the multilinear polytope with hypergraphs,
we investigate the connections between hypergraph acyclicity and the complexity of the facial structure of the multilinear polytope. We characterize the acyclic hypergraphs for which a polynomial-size extended formulation for the multilinear polytope can be constructed in polynomial time.
\end{abstract}

\section{Introduction}
\begin{figure}
    \centering
    \hspace{-1.8cm}
    \hspace{\stretch{1}}
    \includegraphics[scale=.7]{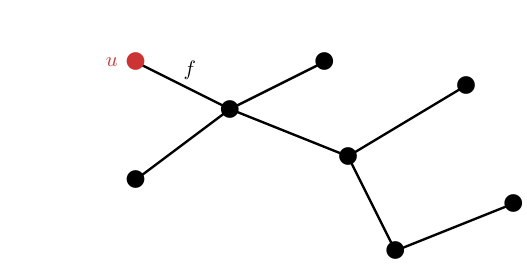}
    \hspace{\stretch{1}}
    \includegraphics[scale=.7]{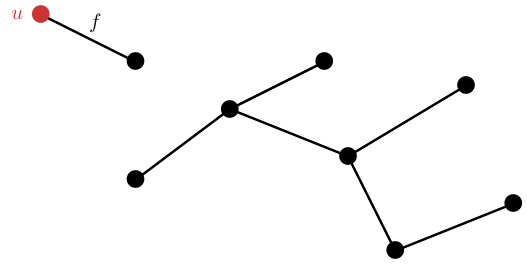}
    \hspace{\stretch{1}}
    \caption{An illustration of the proof technique for~\cref{padberg}.}
    \label{fig tree}
\end{figure}

Binary polynomial optimization, \ie the problem of finding a binary point maximizing a polynomial function, is a fundamental $\NP$-hard problem in discrete optimization with a wide range of applications across science and engineering. To formally define this problem, we employ a hypergraph representation scheme introduced in~\cite{dPKha17MOR}. A \emph{hypergraph} $G$ is a pair $(V,E)$, where $V$ is a finite set of nodes and $E$ is a set of subsets of $V$, called the edges of $G$. Throughout this article, we consider hypergraphs without loops or parallel edges, in which case $E$ is a set of subsets of $V$ of cardinality at least two. Moreover, the \emph{rank} of a hypergraph is the maximum cardinality of any edge in $E$.
With any hypergraph $G= (V,E)$, and costs $c_v$, $v \in V$, and $c_e$, $e \in E$, we associate the following binary polynomial optimization problem:
\begin{align}
\label[problem]{prob BPO}
\tag{BPO}
\begin{split}
\max & \qquad \sum_{v\in V} {c_v z_v} + \sum_{e\in E} {c_e \prod_{v\in e} {z_v}} \\
{\rm s.t.} & \qquad z_v \in \{0,1\} \qquad \forall v \in V,
\end{split}
\end{align}
where without loss of generality we assume $c_e \neq 0$ for all $e \in E$. We then proceed with linearizing the objective function by introducing a new variable for each product term to obtain an equivalent reformulation of \cref{prob BPO} in a lifted space of variables:
\begin{align}
\label[problem]{prob l-BPO}
\tag{$\ell$-BPO}
\begin{split}
\max & \qquad \sum_{v\in V} {c_v z_v} + \sum_{e\in E} {c_e z_e} \\
{\rm s.t.} & \qquad z_e = \prod_{v\in e} {z_v} \qquad \forall e \in E\\
& \qquad z_v \in \{0,1\} \qquad \forall v \in V.
\end{split}
\end{align}
To solve \cref{prob l-BPO} efficiently using polyhedral techniques, it is essential to understand the facial structure of the convex hull of its feasible region.
To this end, in the same vein as~\cite{dPKha17MOR}, we define the~\emph{multilinear set} as
\begin{equation*} \label{eq: SG}
\S(G)= \Big\{ z \in \{0,1\}^{V \cup E} : z_e = \prod_{v \in e} {z_{v}}, \; \forall e \in E \Big\},
\end{equation*}
and we refer to its convex hull as the~\emph{multilinear polytope} and denote it by $\MP(G)$. 
A simple polyhedral relaxation of $\S(G)$ can be obtained by replacing each term $z_e = \prod_{v \in e} {z_{v}}$ by its convex hull over the unit hypercube:
\begin{equation*}
\MP^{\LP}(G) =\Big\{z: 
z_v \leq 1, \forall v \in V; \; 
z_e \geq 0, \; z_e \geq \sum_{v\in e}{z_v}-|e|+1, \forall e \in E; \;
z_e \leq z_v, \forall e \in E, \forall v \in e\Big\}.
\end{equation*}
The above relaxation is often referred to as the \emph{standard linearization} and has been used extensively in the literature~\cite{yc93}. Notice that $\MP^{\LP}(G)$ is defined by at most $|V|+ (r+2) |E|$ inequalities, where $r$ denotes the rank of $G$.
In the special case with $r = 2$; \ie when all product terms in $\S(G)$ are products of two variables, the multilinear polytope coincides with the well-known Boolean quadric polytope ${\rm BQP}(G)$~\cite{Pad89}. Padberg~\cite{Pad89} gives a characterization of graphs for which the Boolean quadric polytope coincides with its standard linearization.
\begin{theorem}[\cite{Pad89}]\label{padberg}
${\rm BQP}(G) = \MP^{\LP}(G)$ if and only if the graph $G$ is acyclic.    
\end{theorem}
We next outline an alternative proof strategy for the ``if'' direction of~\cref{padberg}, which we will later generalize to obtain our characterizations. 
To this end, we introduce the concept of decomposability of multilinear sets, which plays a key role throughout our results. 
Consider hypergraphs $G_1 = (V_1,E_1)$ and $G_2=(V_2,E_2)$ such that $V_1 \cap V_2 \neq \emptyset$.
We denote by $G_1 \cap G_2$ the hypergraph $(V_1 \cap V_2, E_1 \cap E_2)$.
Let $G := (V_1 \cup V_2, E_1 \cup E_2)$.
We say that $\MP(G)$ is \emph{decomposable into} $\MP(G_1)$ and $\MP(G_2)$ if the system comprised a description of $\MP(G_1)$ and a description of $\MP(G_2)$, is a description of $\MP(G)$.
See~\cite{dPKha18MPA} for a detailed study of the decomposability of multilinear sets.
We now outline the proof strategy for~\cref{padberg}, and we refer the reader to~\cref{fig tree} for an illustration.
Let $G=(V,E)$ be an acyclic graph. Take any leaf $u \in V$ and consider the edge $f$ incident to it. 
Let $G_1 = (f,\{f\})$ and let $G_2=(V \setminus \{u\}, E \setminus \{f\})$. Then the graph $G_1 \cap G_2$ contains only one node, and hence by theorem~1 in~\cite{dPKha18MPA}, we deduce that ${\rm BQP}(G)$
is decomposable into ${\rm BQP}(G_1)$ and ${\rm BQP}(G_2)$.
By definition, ${\rm BQP}(G_1)$ is given by its standard linearization. The proof then follows by induction on the number of nodes of $G$.

Hence, it is natural to ask whether the multilinear polytope of acyclic hypergraphs has a simple structure as well.
Interestingly, unlike graphs, the notions of cycles and acyclicity in hypergraphs are not unique. In fact, the notion of graph acyclicity has been extended to several different \emph{degrees of acyclicity} of hypergraphs.
The most well-known types of acyclic hypergraphs, in increasing order of generality, are Berge-acyclic, $\gamma$-acyclic, $\beta$-acyclic, and $\alpha$-acyclic hypergraphs~\cite{fagin83,BeeFagMaiYan83,Dur12}. Next, we define these notions. 

A \emph{Berge-cycle} in $G$ of length $t$ for some $t \ge 2$, is a sequence $C = v_1, e_1, v_2, e_2,\ldots,v_t,e_t,v_1$, where $v_i$, $i\in \{1,\dots,t\}$ are distinct nodes of $G$, $e_i$ , $i\in \{1,\dots,t\}$ are distinct edges of $G$, $v_i, v_{i+1} \in e_i$ for $i \in \{1,\dots, t-1\}$, and $v_1, v_t \in e_t$. A hypergraph is called \emph{Berge-acyclic} if it contains no Berge-cycles. A \emph{$\gamma$-cycle} in $G$ is a Berge-cycle $C = v_1, e_1, v_2, e_2,\dots,v_t,e_t,v_1$ such that $t \ge 3$, and the node $v_i$ belongs to $e_{i-1}$, $e_{i}$ and no other $e_j$, for all $i \in \{2,\dots, t\}$. A hypergraph is called \emph{$\gamma$-acyclic} if it contains no $\gamma$-cycles. A \emph{$\beta$-cycle} in $G$ is a $\gamma$-cycle $C = v_1, e_1, v_2, e_2,\ldots,v_t,e_t,v_1$ such that the node $v_1$ belongs to $e_1$, $e_t$ and no other $e_j$.
A hypergraph is called \emph{$\beta$-acyclic} if it contains no $\beta$-cycles. Finally, let us define $\alpha$-acyclic hypergraphs; while~\cite{JegNdi09} offers a definition of $\alpha$-cycles, this definition is fairly complicated and does not provide any intuitive connection between $\alpha$-cycles and other types of cycle defined above. Hence, we present a characterization for $\alpha$-acyclicity that is considerably simpler to understand and is useful for our results.
A set $F$ of subsets of a finite set $V$ has the \emph{running intersection property} if there exists an ordering $p_1, p_2, \ldots, p_m$ of the sets in $F$ such that
for each $k = 2, \dots,m$, there exists $j < k$ such that $p_k \cap \big(\bigcup_{i < k}{p_i}\big) \subseteq p_j$.
A hypergraph is \emph{$\alpha$-acyclic} if and only if its edge set has the running intersection property~\cite{BeeFagMaiYan83}. 

\medskip


In this article, we characterize the acyclic hypergraphs for which a polynomial-size extended formulation for the multilinear polytope can be constructed in polynomial time~\cite{dPKha18SIOPT,dPKha23mMPA}.

\section{The multilinear polytope of $\alpha$-acyclic hypergraphs}

As we mentioned before, $\alpha$-acyclic hypergraphs are the most general type of acyclic hypergraphs. In fact, given any hypergraph $G$, the hypergraph obtained from $G$ by adding one edge that contains all nodes of $G$ is $\alpha$-acyclic. It is therefore not surprising that solving~\cref{prob BPO} over $\alpha$-acyclic hypergraphs is, in general, $\NP$-hard.

\begin{theorem}[\cite{dPDiG23ALG}]\label{thhard}
\cref{prob BPO} over $\alpha$-acyclic hypergraphs is strongly $\NP$-hard. Furthermore, it is
$\NP$-hard to obtain an $\kappa$-approximation for~\cref{prob BPO}, with $\kappa > \frac{16}{17} \approx 0.94$.
\end{theorem}
\cref{thhard} implies that, unless $\P = \NP$, one cannot construct, in polynomial time, a polynomial-size extended formulation for the multilinear polytope of $\alpha$-acyclic hypergraphs. 
However, by making further assumptions about the rank of $\alpha$-acyclic hypergraphs, one can obtain a polynomial-size extended formulation for the multilinear polytope.
Given a hypergraph $G$, in the following we say that an edge of $G$ is \emph{maximal}, if it is not contained in any other edge of $G$.  Moreover, given a hypergraph $G= (V, E)$ and a subset $V' \subset V$, the \emph{section hypergraph} of $G$ induced by $V'$ is the hypergraph $G'=(V',E')$, where $E'=\{e \in E: e \subseteq V'\}$.
The \emph{reduction} of a hypergraph $G=(V,E)$ is the hypergraph $(V, F)$, where $F$ is the set of maximal edges of $G$.
It is well-known that a hypergraph is $\alpha$-acyclic if and only if its reduction is $\alpha$-acyclic \cite{BeeFagMaiYan83}.

\begin{theorem}\label{thalpha}
Let $G = (V, E)$ be an $\alpha$-acyclic hypergraph of rank $r$. Then $\MP(G)$ has an extended formulation with at most $2^r \min\{|V|,|F|\}$ variables and inequalities, where $F$ is the set of maximal edges of $G$.
Moreover, all coefficients and right-hand side constants in the extended formulation are $0,\pm 1$.
\end{theorem}

By~\cref{thalpha}, if $G$ is an $\alpha$-acyclic hypergraph of rank $r$, with $r = O(\log \poly(|V|, |E|))$, then $\MP(G)$ has a polynomial-size extended formulation, where by $\poly(|V|, |E|)$, we imply a polynomial function in $|V|,|E|$.

\begin{figure}
    \centering
    \hspace{-1.2cm}
    \hspace{\stretch{1}}
    \includegraphics[scale=.55]{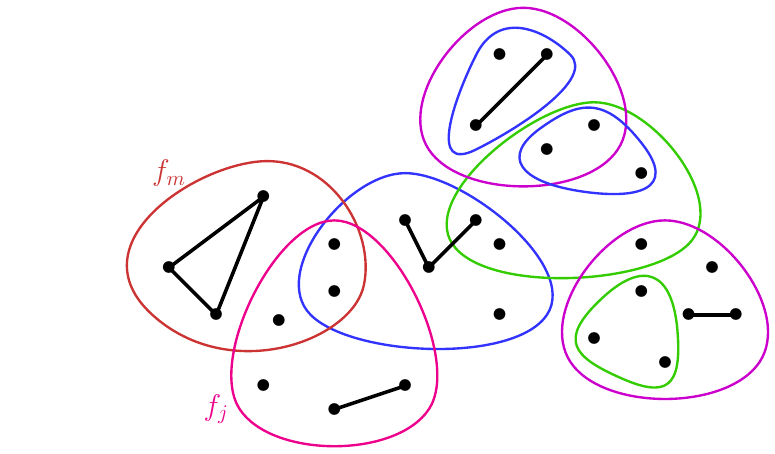}
    \hspace{\stretch{3}}
    \includegraphics[scale=.55]{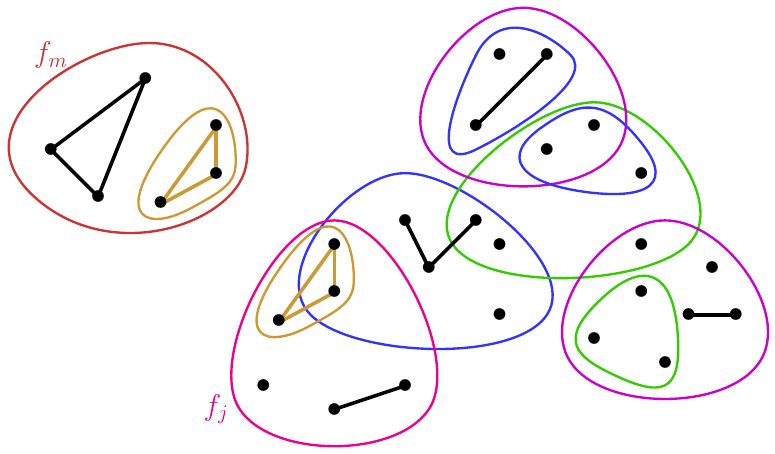}
    \hspace{\stretch{1}}
    \caption{An illustration of the proof technique for~\cref{thalpha}.}
    \label{fig alpha}
\end{figure}

Next, we provide a proof sketch of~\cref{thalpha}.
See~\cref{fig alpha} for an illustration of this proof technique. 
%
Let $R=(V,F)$ be the reduction of $G$, which is also $\alpha$-acyclic.
Then $F$ has the running intersection property.
Let $f_1, f_2, \ldots, f_m$ be a running intersection ordering of $F$.
In particular, there exists $j < m$, such that $f_m \cap\left(f_1 \cup f_2 \cup \cdots \cup f_{m-1}\right) \subseteq f_j$.
Let $R_2$ be the section hypergraph of $R$ induced by $f_1 \cup f_2 \cup \cdots \cup f_{m-1}$.
We can show that $R_2$ is $\alpha$-acyclic and it does not contain $f_m$.
Denote by $H$ the hypergraph obtained from $G$ by adding all edges contained in $f_j \cap f_m$.
To obtain an extended formulation of $\MP(G)$, it suffices to derive an extended formulation of $\MP(H)$.
Let $H_1$ be the section hypergraph of $H$ induced by $f_m$ and let $H_2$ be the section hypergraph of $H$ induced by $f_1 \cup f_2 \cup \cdots \cup f_{m-1}$.
We observe that $H_1 \cup H_2 = H$ and that $H_1 \cap H_2$ is a complete hypergraph.
Hence, by theorem~1 in~\cite{dPKha18MPA}, $\MP(H)$ is decomposable into $\MP(H_1)$ and $\MP(H_2)$.
The hypergraph $H_1$ has at most $r$ nodes. 
Therefore, $\MP(H_1)$ has an extended formulation with at most $2^r$ variables and inequalities~\cite{SheAda90}.
Moreover, $H_2$ is an $\alpha$-acyclic hypergraph of rank at most $r$. 
By repeating the above argument and noticing that 
$H_2$ hast at most $|V|-1$ nodes, and $m-1$ maximal edges, we deduce that the resulting extended formulation has at most $2^r\min\{|V|,|F|\}$ variables and inequalities.

In~\cite{WaiJor04,Lau09,BieMun18}, the authors give extended formulations for the convex hull of the feasible set of (possibly constrained) binary polynomial optimization problems. The size of these extended formulations is parameterized in terms of the ``treewidth'' of their so-called intersection graphs. See for example~\cite{BieMun18} for definitions of treewidth and intersection graphs.
For the unconstrained case, their result can be stated as follows:
\begin{theorem}\label{thtw}
Let $G=(V,E)$ be a hypergraph and let $w$ be the treewidth of its intersection graph. Then $\MP(G)$ has an extended formulation
with at most $2^w |V|$ variables and inequalities.
\end{theorem}
As we detail in~\cite{dPKha21MOR}, \cref{thalpha} is equivalent to~\cref{thtw}. 

We next show that even for a quadratic objective function, a bounded treewidth is a necessary condition for tractability of~\cref{prob BPO}.
To this end, we first introduce some terminology.
Given an instance $\Lambda$ of \cref{prob BPO}, we refer to the \emph{size} of $\Lambda$ (also known as \emph{bit size}, or \emph{length}), denoted by $\|\Lambda\|$, as the number of bits required to encode it. We say that a countable family of graphs $\{G_k\}_{k=1}^\infty$ is \emph{polynomial-time enumerable} if a description of $G_k$ is computable in time bounded by a polynomial in $k$. 
Recall that the Bounded-error Probabilistic Polynomial time ($\BPP$) is the class of decision problems solvable by a probabilistic Turing machine in polynomial time with an error probability bounded by $1/3$ for all instances. 
Our intractability result is under the assumption that $\NP \not\subseteq \BPP$.
It is widely believed that $\P=\BPP$, which implies that $\NP \not\subseteq \BPP$ is equivalent to $\P \neq \NP$.
For a precise definition of $\BPP$ and the commonly believed $\NP \not\subseteq \BPP$ hypothesis, we refer the reader to \cite{AroBar09}.
In~\cite{dPKha24limits}, we obtain a necessary condition for the tractability of~\cref{prob BPO}:


\begin{theorem}
\label{th treewidth}
Let $\{G_k\}_{k=1}^\infty$ be a polynomial-time enumerable family of graphs where the treewidth of $G_k$ equals $k$, for all $k$.
Let $f$ be an algorithm that solves any instance $\Lambda_k$ of~\cref{prob BPO} on graph $G_k$ in time $T(k) \cdot \poly(\|\Lambda_k\|)$.
Then, assuming $\NP \not\subseteq \BPP$, $T(k)$ grows super-polynomially in $k$. 
\end{theorem}

Due to the connection between the rank of $\alpha$-acyclic hypergraphs and the treewidth of intersection graphs~\cite{dPKha21MOR}, \cref{th treewidth} implies that for general $\alpha$-acyclic hypergraphs a bounded rank is a necessary condition for tractability of~\cref{prob BPO}. This in turn implies that, unless $\NP = \BPP$, one cannot obtain in polynomial time, a polynomial-size extended formulation for the multilinear polytope of any sequence of $\alpha$-acyclic hypergraphs of increasing rank.

\section{Polynomial-size extended formulations}
While for $\alpha$-acyclic hypergraphs, a polynomial-size extended formulation for $\MP(G)$ can be obtained only if the rank of the hypergraph is ``bounded'',  as we detail in this section, for the remaining types of acyclic hypergraph no such restrictive assumptions are required. 
Recall that acyclic hypergraphs in increasing degree of generality are:
$$
\text{Berge-acyclic} \subset \text{$\gamma$-acyclic} \subset \text{$\beta$-acyclic} \subset \text{$\alpha$-acyclic}.
$$
In~\cite{dPKha23mMPA} we obtain a polynomial-size extended formulation for the multilinear polytope of $\beta$-acyclic hypergraphs:

\begin{theorem}\label{thbeta}
Let $G = (V,E)$ be a $\beta$-acyclic hypergraph of rank $r$.
Then there exists a polynomial-size extended formulation of $\MP(G)$ comprising at most
$(3r-4)|V|+4|E|$ inequalities and at most $(r-1) |V|+|E|$ variables.
Moreover, all coefficients and right-hand side constants in the extended formulation are $0,\pm 1$.
\end{theorem}

The proof of~\cref{thbeta} relies on the key concept of nest points of hypergraphs. A node $v \in V$ is a \emph{nest point} of $G$ if the set of edges of $G$ containing $v$ is totally ordered. We define the hypergraph obtained from $G = (V,E)$ by \emph{removing} a node $v \in V$ as $G - v := (V',E')$, where $V' := V \setminus \{ v \}$ and $E' := \{ e \setminus \{v\} : e \in E, \ |e \setminus \{v\}| \ge 2 \}$.
A \emph{nest point elimination order} is an ordering $v_1, \dots, v_n$ of the nodes of $G$, such that $v_1$ is a nest point of $G$, $v_2$ is a nest point of $G - v_1$, and so on, until $v_n$ is a nest point of $G - v_1 - \dots - v_{n-1}$. 
To prove~\cref{thbeta}, we use the following characterization of $\beta$-acyclic hypergraphs: a hypergraph $G=(V, E)$ is $\beta$-acyclic if and only if it has a nest point elimination order~\cite{Dur12}. 

\begin{figure}
    \centering
    \hspace{-2.4cm}
    \hspace{\stretch{1}}
    \includegraphics[scale=.55]{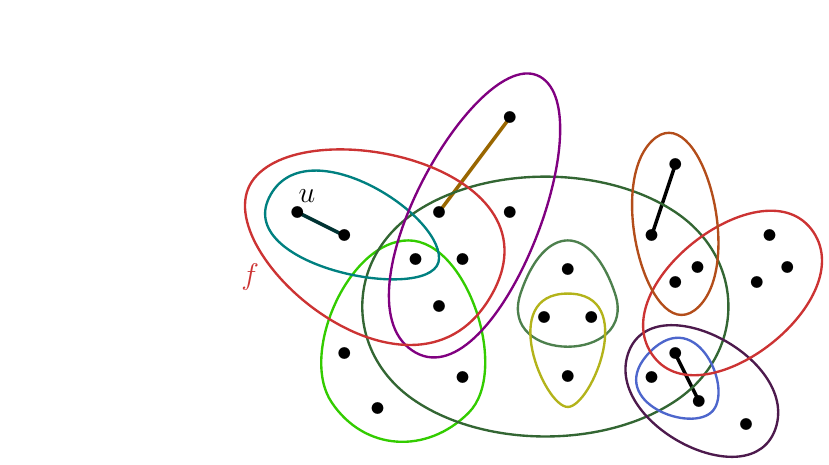}
    \hspace{\stretch{3}}
    \includegraphics[scale=.55]{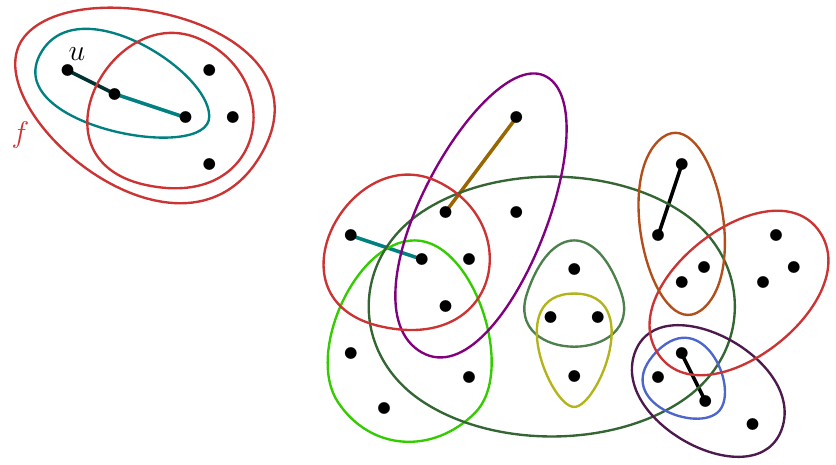}
    \hspace{\stretch{1}}
    \caption{An illustration of the proof technique for~\cref{thbeta}.}
    \label{fig beta}
\end{figure}

Next, we outline the proof strategy for~\cref{thbeta}. See~\cref{fig beta} for an illustration of this proof technique.  
Let $G=(V, E)$ be a $\beta$-acyclic hypergraph and let $u \in V$ be a nest point of $G$. Let $F$ denote the set of edges containing $u$ and denote by $f$ the edge of the largest cardinality in $F$. Define $P := \{e \setminus \{u\}: e \in F,  |e| \geq 3\}$, and for each $p \in P$, if $p$ is not an edge of $G$, then add it to $G$ and call the resulting hypergraph $G'$. Then by theorem~4 in~\cite{dPKha23mMPA}, the multilinear polytope $\MP(G')$ is decomposable into $\MP(\bar G)$ and $\MP(G-u)$, where $\bar G = (f, F \cup P)$. Clearly, the hypergraph $\bar G$ has a very special structure. In~\cite{dPKha23mMPA}, we prove that $\MP(\bar G)$ is defined by $5|f|+2$ inequalities in the original space of variables. Moreover, it can be shown that $\MP(G-u)$ is a $\beta$-acyclic hypergraph. Therefore, the proof follows by induction on the number of nodes of $G$. 

It is important to note that the standard linearization $\MP^{\LP}(G)$ often leads to very weak relaxations of $\MP(G)$ for $\beta$-acyclic hypergraphs. \cref{thbeta} implies that while the proposed extended formulation for $\MP(G)$ contains $(r-2)|V|$ additional variables, it has fewer inequalities than the standard linearization if $|E| \geq 3 |V|$. 
We should also remark that the inequalities defining our proposed extended formulation are very sparse; that is, they
contain at most four variables with non-zero coefficients; a feature that is highly beneficial from a computational perspective.

\section{Explicit characterizations in the original space}
In the previous section, we presented a polynomial-size extended formulation for the multilinear polytope of $\beta$-acyclic hypergraphs. It is often desirable to obtain an explicit description of the multilinear polytope in the original space. 
In principle, one could project out the extra variables from the extended formulation.
In this section, we review our results regarding the structure of the multilinear polytope in the original space for Berge-acyclic, $\gamma$-acyclic and $\beta$-acyclic hypergraphs~\cite{dPKha18SIOPT,dPKha23mMPA}. As an important byproduct of these results, we obtain new classes of cutting planes that enable us to construct stronger linear programming relaxations for general mixed-integer polynomial optimization problems~\cite{dPKha21MOR,dPKhaSah20MPC}.

By~\cref{padberg}, the Boolean quadric polytope coincides with its standard linearization if and only if $G$ is an acyclic graph. In~\cite{dPKha18SIOPT}, we obtain a generalization of this result:

\begin{theorem}\label{thberge}
    $\MP(G) = \MP^{\LP}(G)$ if and only if $G$ is a Berge-acyclic hypergraph.
\end{theorem}

Next, we define a class of valid inequalities for the multilinear set, which enables us to characterize the multilinear polytope of $\gamma$-acyclic hypergraphs in the original space. Let $G=(V,E)$ be a hypergraph, let $e_0$ be an edge of $G$, and denote by $e_k$, $k \in K$, the set of all edges of $G$ adjacent to $e_0$.
Let $T$ be a nonempty subset of $K$ such that
\begin{equation}\label{cond}
\Big|(e_0 \cap e_i)\setminus \bigcup_{j \in T \setminus \{i\}}{(e_0 \cap e_j)}\Big| \geq 2, \quad \forall i \in T.
\end{equation}
Then the~\emph{flower inequality} centered at $e_0$ with neighbors $e_k$, $k \in T$, is given by:
\begin{equation}\label{flowerIneq}
\sum_{v \in e_0\setminus \cup_{k\in T} {e_k}}{z_v}+\sum_{k \in T} {z_{e_k}} - z_{e_0} \leq |e_0\setminus \cup_{k\in T} {e_k}|+|T|-1.
\end{equation}
We refer to inequalities of the form~\eqref{flowerIneq}, for all nonempty $T \subseteq K$
satisfying condition~\eqref{cond}, as the system of flower inequalities centered at $e_0$.
We then define the \emph{flower relaxation} $\MP^{F}(G)$ as the polytope obtained by adding the system of flower inequalities centered at each edge of $G$ to $\MP^{\LP}(G)$. In~\cite{dPKha18SIOPT} we examine the tightness of the flower relaxation:
\begin{theorem}\label{thgamma}
    $\MP(G) = \MP^{F}(G)$ if and only if $G$ is a $\gamma$-acyclic hypergraph.
\end{theorem}
In~\cite{dPKha18SIOPT} we prove that the multilinear polytope of $\gamma$-acyclic hypergraphs may contain exponentially many facet-defining inequalities, in general. Nonetheless, we devise a strongly polynomial-time algorithm to separate flower inequalities over the multilinear polytope of $\gamma$-acyclic hypergraphs. In~\cite{dPKhaSah20MPC} we prove that the separation problem for flower inequalities over general hypergraphs is $\NP$-hard. 
However, for hypergraphs with fixed rank, the separation problem can be solved in $O(|E|^2)$ operations. In~\cite{dPKha21MOR}, we propose a new class of cutting planes, which we refer to as~\emph{running intersection inequalities}. Under certain assumptions, these inequalities are stronger than flower inequalities. In~\cite{dPKhaSah20MPC}, we incorporate flower inequalities and running intersection inequalities in the state-of-the-art mixed-integer nonlinear programming solver {\tt BARON}~\cite{IdaNick18}. Results show that the proposed cutting planes significantly improve the performance of {\tt BARON} for both random test sets and image restoration problems.

Finally, let us discuss the facial structure of the multilinear polytope of $\beta$-acyclic hypergraphs.
As $\beta$-acyclicity subsumes $\gamma$-acyclicity, it follows that also the multilinear polytope of a $\beta$-acyclic hypergraph may contain exponentially many facet-defining inequalities.
However, as we detail in~\cite{dPKha23mMPA}, unlike the multilinear polytope of Berge-acyclic and $\gamma$-acyclic hypergraphs, the multilinear polytope of $\beta$-acyclic hypergraphs may contain \emph{very dense} facets, in general. That is, inequalities containing as many as $\theta(|E|)$ nonzero coefficients even when $r=8$. From a computational perspective, sparsity is key to the effectiveness of cutting planes in a branch-and-cut framework. Indeed, all existing families of cutting planes for multilinear sets, such as flower inequalities~\cite{dPKha18SIOPT} and running intersection inequalities~\cite{dPKha21MOR} are sparse. Namely, for a rank $r$ hypergraph, flower inequalities contain at most $\frac{r}{2}$ nonzero coefficients, and running intersection inequalities contain at most $2(r-1)$ nonzero coefficients.  This is significant since, for most multilinear sets appearing in applications, we have $r \ll |E |$. To conclude, although obtaining an explicit description for the multilinear polytope of $\beta$-acyclic hypergraphs in the original space remains an open question, we suspect that due to the issues described above, such a characterization may not have any practical impact.

\paragraph{Acknowledgments.}
We would like to thank the INFORMS Computing Society (ICS) Prize 
Committee consisting of Mirjam~Dur, Yufeng~Liu, Jim~Luedtke, and Uday~Shanbhag (chair). It is a great honor to have been selected as recipients of this award. A.~Del~Pia was supported in part by AFOSR grant FA9550-23-1-0433, and A.~Khajavirad was supported in part by AFOSR grant FA9550-23-1-0123.

\begin{footnotesize}
\bibliographystyle{plain}

\end{footnotesize}
\end{document}